\documentclass{amsart}

\usepackage{amsfonts}
\usepackage{amsmath}
\usepackage{amsthm}
\usepackage{amssymb}
\usepackage[english]{babel}
\usepackage{graphics}
\usepackage{latexsym}
\usepackage{longtable} \usepackage{mathrsfs}
\usepackage{graphicx}
\usepackage{wrapfig} 
\usepackage[pdftex,colorlinks=true]{hyperref}

\newtheorem{theorem}{Theorem}
\newtheorem{corollary}[theorem]{Corollary}
\newtheorem{lemma}[theorem]{Lemma}

\newtheorem{claim}[theorem]{Claim}
\newtheorem{example}[theorem]{Example}

\theoremstyle{definition}

\newcommand{\ess}{\textrm{ess}}

\newcommand{\spn}{\textrm{span}}

\newcommand{\BPL}{\medskip \noindent \textbf{Proof of Lemma} }

\newcommand{\BPT}{\medskip \noindent \textbf{Proof of Theorem} }

\newcommand{\sub}{\subseteq}

\newcommand{\noi}{\noindent}
\newcommand{\ms}{\medskip}

\newcommand{\Om}{\Omega}

\newcommand{\la}{\lambda}

\newcommand{\e}{\varepsilon}
\newcommand{\de}{\delta}
\newcommand{\al}{\alpha}
 \newcommand{\larrow}{\longrightarrow}

\newcommand{\R}{\mathbb{R}}

\newcommand{\bt}{\begin{theorem}}\newcommand{\et}{\end{theorem}}
\newcommand{\bl}{\begin{lemma}}\newcommand{\el}{\end{lemma}}
\newcommand{\be}{\begin{equation}}\newcommand{\ee}{\end{equation}}
\newcommand{\bc}{\begin{claim}}\newcommand{\ec}{\end{claim}}
\newcommand{\bp}{\begin{proof}}\newcommand{\ep}{\end{proof}}

\begin{document}

\title{Explicit Singular Viscosity Solutions of the Aronsson Equation}

\author{\textsl{Nikolaos I. Katzourakis}}
\address{BCAM - Basque Center for Applied Mathematics, Biskaia Technology Park, Building 500, E-48160, Derio, Spain}
\email{nkatzourakis@bcamath.org}





\begin{abstract}
\noi We establish that when $n\geq 2$ and $H \in C^1(\R^n)$ is
a Hamiltonian such that some level set contains a line segment, the
Aronsson equation $D^2  u : H_p(Du) \otimes H_p(Du)= 0$ admits explicit entire viscosity solutions. They are superpositions
of a linear part plus a Lipschitz continuous singular part which in general is non-$C^1$ and nowhere twice differentiable. In particular, we supplement the $C^1$ regularity result of Wang and Yu \cite{W-Y} by deducing that strict level convexity is necessary for $C^1$ regularity of solutions.
\ms

\noi \textsc{R\'esum\'e}. Nous d\'emontrons que, pour $n \geq 2$ et un Hamiltonien $H\in {C}^1(\mathbb{R}^n)$ tel qu'au moins une de ses lignes de niveau contienne un segment de droite, l'\'equation de Aronsson $D^2  u : H_p(Du) \otimes H_p(Du)= 0$ admet des solutions de viscosit\'e explicites d\'efinies sur $\mathbb{R}^n$. Elles sont superpositions d'une partie lin\'eaire et d'une partie continue, lipschitzienne, singuli\`ere qui, en g\'en\'eral, n'est pas $C^1$ et est nulle part deux fois d\'erivable. Plus pr\'ecis\'ement, nous compl\'etons le r\'esultat de r\'egularit\'e \'etablit par Wang et Yu \cite{W-Y} en montrant que la stricte convexit\'e des lignes de niveau est n\'ecessaire pour que les solutions soient $C^1$.
\end{abstract}

\maketitle

\section{Introduction}

\noi Let $H \in C^1(\R^n)$ be a Hamiltonian function and $n\geq 2$.
We discuss aspects of the $C^1$ regularity problem of viscosity
solutions to the Aronsson PDE, which is defined on $u
\in C^2(\R^n)$ by
 \begin{align}\label{eq1}
\mathcal{A}[u]\ :=&\ D^2  u : H_p(Du) \otimes H_p(Du)\ = \ 0.
 \end{align}
Here, $\mathcal{A}[u]$ means $\sum_{i,j=1}^n D^2_{ij}u \,
H_{p_i}(Du)\, H_{p_j}(Du)$ and $H_{p_i}=D_{p_i}H$. \eqref{eq1} defines a quasilinear highly degenerate elliptic PDE. It
arises in $L^\infty$ variational problems of the supremal functional
$E_\infty(u,\Om) := \|H (Du)\|_{L^{\infty}(\Om)}$, as well as in
other contexts (Barron-Evans-Jensen \cite{BEJ}). When
$H(p)=\frac{1}{2}|p|^2$, \eqref{eq1} reduces to the
$\infty$-Laplacian  $\Delta_{\infty}u := D^2  u : Du \otimes Du = 0$.
Under reasonable convexity, coercivity and regularity assumptions on
$H$, there exists a unique continuous solution of the Dirichlet
problem with Lipschitz boundary data, interpreted in the viscosity
sense of Crandall-Ishii-Lions \cite{CIL} which actually is Lipschitz continuous.
However, the $C^1$ regularity problem for \eqref{eq1} remains open.
Wang and Yu \cite{W-Y} established that when $n=2$, $H$ is in
$C^2(\R^2)$ with $H \geq H(0)=0$ and it is uniformly convex, then viscosity solutions of \eqref{eq1} over $\Om \sub \R^2$
are in $C^1(\Om)$. When $n> 2$, viscosity solutions are linearly
approximatable in the sense of De Pauw-Koeller
\cite{DePK}, having approximate gradients.

Herein we prove that when a level set $\{H=c\}$ of $H$ contains a
straight line segment, there exists an entire viscosity solution of
\eqref{eq1} given as superposition of a linear term plus a rather
arbitrary Lipschitz continuous term. The latter may \emph{not} be $C^1$; moreover, it may well be nowhere twice differentiable with Hessian realized only as a singular distribution and not even as a Radon measure.

We note that our \emph{only} assumption is $H$ being constant along
a line segment but arbitrary otherwise. This suffices for these
solutions to appear. Actually, they arise as a.e. solutions of the Hamilton-Jacobi equation $H(Du) = c$. However, we work with the second order PDE \eqref{eq1} ignoring the relation between viscosity solutions of \eqref{eq1} and solutions of $H(Du) = c$. We just notice that in the $C^2$ context, the
identity
 \be \label{eq4}
D^2  u : H_p(Du) \otimes H_p(Du) \ = \ H_p(Du)\cdot D\big(H(Du)\big)
 \ee
suffices to imply $\mathcal{A}[u]=0$, whenever $H(Du)=c$. Let us now
state our result.

\begin{theorem}\label{th1} We assume that $H \in
C^1(\R^n)$, $n\geq 2$ and there exists a straight line segment
$[a,b] \sub \R^n$ along which $H$ is constant. Then, for any $f \in W^{1,\infty}_{loc}(\R)$ satisfying $\|f'\|_{L^\infty(\R)}<1$, the formula
 \be \label{eq5}
 u(x) \ :=\ \frac{b+a}{2} \cdot x \ +
     \ f\left(\frac{b-a}{2} \cdot x\right),\ \ \
     x \in \R^n,
 \ee
defines an entire viscosity solution $u \in W^{1,\infty}_{loc}(\R^n)$ of the
Aronsson equation.
\end{theorem}
\noi We deduce that the existence of the non-$C^1$ solutions
\eqref{eq5} implies the following
\begin{corollary}
Strict level convexity of the Hamiltonian $H$ is necessary to obtain
$C^1$ regularity of viscosity solutions to the Aronsson PDE in all
dimensions $n\geq 2$.
\end{corollary}
In particular, the uniform convexity assumption of Wang and Yu
\cite{W-Y} can not be relaxed to mere convexity, unless if strict
level-convexity is additionally assumed.

We observe that $C^1$ regularity of solutions is not an issue of
regularity of $H$; the singular solutions \eqref{eq5} persist even
when $H\in C^\infty(\R^n)$. The sensitive dependence of regularity
on the convexity of $H$ is a result of the geometric degeneracy
structure of the PDE $\mathcal{A}[u]= 0$ which in view of \eqref{eq4} can be rewritten as the perpendicularity condition $H_p(Du)\ \bot\ D\big(H(Du)\big)$. Also, the singular solutions persist for arbitrarily small straight line segments, as long as the segments do not trivialize to a point.

\section{Proofs}

\noi For the definition and the properties of viscosity solutions we refer to Crandall-Ishii-Lions \cite{CIL}. We will first prove Theorem \eqref{th1} for smooth functions $f$ and then deduce the full result by approximation.

\begin{lemma}\label{le1} Let $u$ be given by \eqref{eq5} with $f\in C^2(\R)$ satisfying $\|f'\|_{L^\infty(\R)}<1$. Then,

\noi (i) $Du(\R^n) \sub (a,b)$, i.e.\ the range
 of its gradient $Du$ is valued in the open segment $(a,b)=\{x\in \R^n \; |\; x=\la a +(1-\la)b,\; \la\in(0,1)\}$,

\noi (ii) $H_p\big(Du(\R^n)\big)\sub (\spn[b-a])^\bot$, i.e.\ the gradient of $H$ restricted on $Du(\R^n)$ is normal to $(a,b)$.
\end{lemma}

\BPL \ref{le1}. By differentiating \eqref{eq5}, we have
 \be \label{eq11}
Du(x)\ = \ \frac{b+a}{2} \ + \ \frac{1}{2}\
f'\left(\frac{b-a}{2} \cdot x\right)(b-a),
 \ee
 for all $x\in \R^n$. By rearranging \eqref{eq11}, we have
\be
Du(x) =  \left(\frac{1}{2} - \frac{1}{2}
f'\left(\frac{b-a}{2}\cdot x \right)\right) a\ +\ \left[1-  \left(\frac{1}{2} - \frac{1}{2}
f'\left(\frac{b-a}{2}\cdot x \right)\right)\right]b.
 \ee
Since $\|f'\|_{L^\infty(\R)}<1$, there exists a $\de>0$ such that
\be
\de\ \leq \ \frac{1}{2} - \frac{1}{2}
f'\left(\frac{b-a}{2}\cdot x \right) \ \leq \ 1- \de,
\ee
 for all $x\in \R^n$. Hence, $Du(x)$ is for all
$x\in \R^n$ a strict convex combination of $a$ and $b$. Thus, (i) follows. Since $H$ is constant on $[a,b]$, there exists $c\in \R$ such that, for all $t\in(0,1)$, we have 
\be
H\big(t b+(1-t)a\big)\ =\ c.
\ee
Since $H \in C^1(\R^n)$, we may differentiate to find
 \be \label{eq14}
\frac{d}{dt}\Big(H\big(t b+(1-t)a\big)\Big) \ = \ (b-a)\cdot
H_p\big(t b+(1-t)a\big),
 \ee
for all for $0<t<1$. Hence, we obtain that 
$(b-a)\cdot H_p(q) = 0$ for all $q \in (a,b)$. Since by (i) we have $Du(\R^n) \sub (a,b)$, (ii) follows as well. \qed

\[
\includegraphics[scale=0.17]{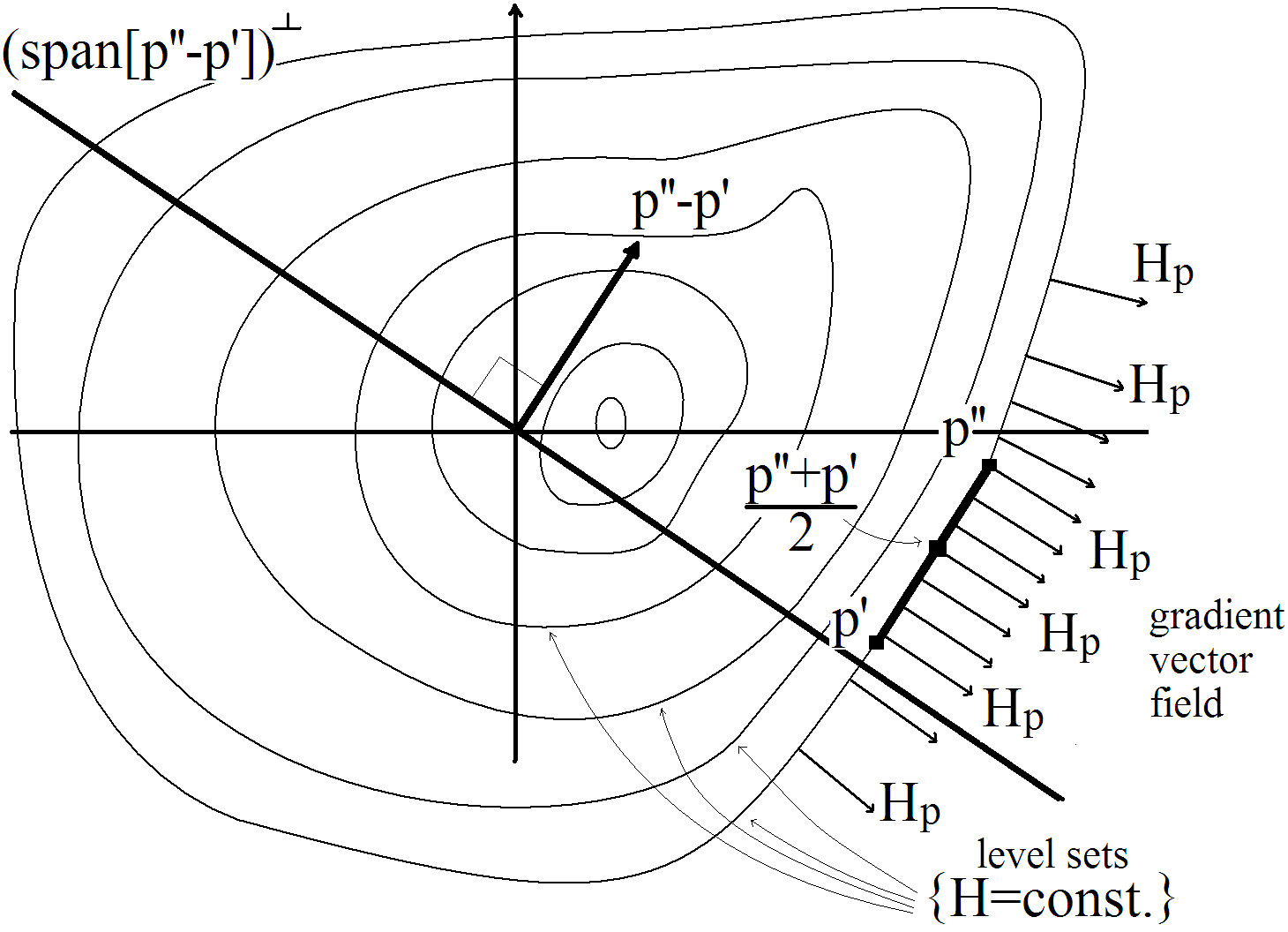}
\]

\begin{lemma} \label{le2} Let $u$ be given by \eqref{eq5} with $f\in C^2(\R)$ satisfying $\|f'\|_{L^\infty(\R)}<1$. Then,  \eqref{eq5}
defines a $C^2(\R^n)$ solution of the Aronsson PDE
\eqref{eq1}.
\end{lemma}

\BPL \ref{le2}. By \eqref{eq5} and our assumption, the Hessian
$D^2u(x)$ exists for all $x\in \R^n$. By differentiating
\eqref{eq11}, we have
 \be \label{eq16}
D^2u(x)\ = \ \frac{1}{4}\
f''\left(\frac{b-a}{2} \cdot x\right)(b-a) \otimes
(b-a).
 \ee
We now calculate using \eqref{eq16} and \eqref{eq11}:
 \begin{align}
 \mathcal{A}[u](x)\ =\ &  D^2u(x):H_p\big(Du(x)\big)\otimes H_p\big(Du(x)\big)\nonumber\\
=\ & \frac{1}{4}\ f''\left(\frac{b-a}{2} \cdot
x\right)(b-a) \otimes (b-a): \nonumber \\
&: H_p\left(\frac{b+a}{2} \ + \ \frac{1}{2}\
f'\left(\frac{b-a}{2} \cdot x\right)(b-a)\right)\otimes \\
&\otimes H_p\left(\frac{b+a}{2} \ + \ \frac{1}{2}\
f'\left(\frac{b-a}{2} \cdot
x\right)(b-a)\right).\nonumber
\end{align}
By employing Lemma \ref{le1}, we have
\begin{align}
 \mathcal{A}[u](x)\ =\ &   \left\{(b-a)\cdot H_p\left(\frac{b+a}{2} +
\frac{1}{2}\, f'\left(\frac{b-a}{2} \cdot
x\right)(b-a)\right)\right\}^2 \cdot \nonumber\\
& \cdot \frac{1}{4}\, f''\left(\frac{b-a}{2} \cdot
x\right) \\
=\ & 0\nonumber
\end{align}
and the Lemma follows. \qed

\ms

\noi Hence, in the case of smooth $u$ the PDE \eqref{eq1} is satisfied because the Hessian $D^2u$ is normal to $H_p(Du)\otimes H_p(Du)$ in the space of symmetric matrices. Now we conclude with the general case of merely Lipschitz $f$. 

\BPT \ref{th1}.  Let $u$ be given by \eqref{eq5} with $f \in W^{1,\infty}_{loc}(\R)$ and $\|f'\|_{L^\infty(\R)}<1$. Let $\eta^\e$, $\e>0$, be the standard mollifier and define $f^\e:=f*\eta^\e \in C^\infty(\R)$. Let also $u^\e$ be given by \eqref{eq5} with $f^\e$ in the place of $f$. Then, $f^\e \larrow f$ in $C^0(\R)$ as $\e\rightarrow 0$ and hence $u^\e \larrow u$ in $C^0(\R^n)$ as $\e\rightarrow 0$. Moreover,
\be
\| {f^{\e}}' \|_{L^\infty(\R)}\ \leq\ \underset{x\in \R}{\ess \, \sup} \int_\R |f'(x-y)|\, |\eta^\e(y)|dy
\ee
and hence $\|{f^{\e}}'\|_{L^\infty(\R)} \leq \|f'\|_{L^\infty(\R)}<1$. Consequently, by Lemmas \ref{le1} and \ref{le2}, all $u^\e$ are smooth entire solutions to the PDE \eqref{eq1}: $\mathcal{A}[u^\e]=0$. By the stability of viscosity solutions, we have $\mathcal{A}[u]=0$ on $\R^n$ in the viscosity sense and Theorem \ref{th1} follows. \qed

\begin{example} The choice $f(t):=\frac{1}{2}|t|$ for $|t|\leq 1$ and $f(t+2)=f(t)$ gives a non-$C^1$ solution $u$ to the PDE \eqref{eq1}. The choice $f(t):=\frac{1}{2}\int_0^t K_{\al,\nu}(s)ds$ with $K_{\al,\nu} \in C^0(\R)$ the singular function of \cite{K} gives a nowhere twice differentiable solution $u$ to \eqref{eq1} with $D^2 u$  existing only as a singular first order distribution.
\end{example}

\medskip

\noindent {\bf Acknowledgement.} {Part of this work was carried out when the Author was a doctoral student at the Department of Mathematics, University of Athens, Greece.}

\bibliographystyle{amsplain}

\end{document}